\documentclass{amsart}
\usepackage{amssymb,amsthm,amsmath,amscd,hyperref}
\def\be{\begin{equation}}
\def\ee{\end{equation}}
\def\ba{\begin{aligned}}
\def\ea{\end{aligned}}
\def\la{\label}

\def\SD{{\mathbb{S}}^{d-1}}
\def\RD{{\mathbb{R}}^{d}}
\def\F{{\mathcal{F}}} 

\def\eps{\epsilon}

\def\ra{\rightarrow}
\def\tr{\qopname \relax o{tr}}

\def\R{{\mathbb{R}}}
\def\N{{\mathbb{N}}}
\def\Z{{\mathbb{Z}}}

\def\min{\qopname \relax o{min}}

\def\Tr{\qopname \relax o{Tr}}
\def\tr{\qopname \relax o{Tr}}

\def\meas{\textrm{mes}}

\def\lp{\left(}
\def\rp{\right)}

\def\HS{{\mathfrak{S}_2}}

\def\spec{\qopname \relax o{spec}}
\def\op{\qopname \relax o{op}}
\def\sgn{\qopname \relax o{sgn}}
\numberwithin{equation}{section}

\newtheorem{theorem}{Theorem}[section]
\newtheorem{lemma}[theorem]{Lemma}
\newtheorem{corollary}[theorem]{Corollary}

\theoremstyle{remark}
\newtheorem{remark}{Remark}[section]

\begin{document}

\title[Szeg\"o limit theorem for operators with discontinuous symbols]
{Szeg\"o limit theorem for operators with discontinuous symbols
and applications to entanglement entropy}
\author{Dimitri Gioev}

\subjclass[2000]{Primary: 81P15, 58J37, 47B35. Secondary: 35S05, 82B10.} 

\address{Department of Mathematics, University of Rochester, 
Hylan Building, Rochester, NY 14627}
\email{gioev@math.rochester.edu}
\begin{abstract}
{
The main result in this paper
is a one term Szeg\"o type asymptotic formula
with a sharp remainder estimate for a class of integral operators
of the pseudodifferential type
with symbols which are allowed to be 
non-smooth or
discontinuous in both position and momentum.
The simplest example of such symbol is the product of the characteristic
functions of two compact sets, one in real space and the other
in momentum space.
The results of this paper are used in 
a study of the violation of the area entropy 
law for free fermions in \cite{GK}.
This work also provides evidence towards a conjecture due to Harold Widom.
}
\end{abstract}
\maketitle
\section{Introduction} 

The problem of estimating entanglement entropy (EE) is currently of 
considerable interest
in the physics community, in particular in condensed matter physics and
in the theory of quantum
information. The interest in EE in condensed matter systems
 is due, in particular, to its scaling behavior and universal properties 
near quantum phase transitions, 
see for example \cite{Cardy,Korepin,Osborne,Osterloh}
and in particular \cite{Vidal}.
EE is an accepted measure of entanglement:
in the quantum information context, entanglement
is necessary for performing quantum 
computations, see 
e.g.~\cite{BBPS,Wootters,Linden,Bennett,Nielsenbook,Vidal1}.
An experimental demonstration of entanglement effects
in a macroscopic system was reported in \cite{Ghosh}.
A connection between EE for spin chain
models and Random Matrix Theory was found in \cite{KM1,KM2}.
 
The evaluation of EE in physical systems of interest
presents considerable
mathematical difficulties \cite{Fannes,JinKorepin,Its,GK}.
In several models EE turns out to be closely related 
to various versions of the strong Szeg\"o
limit theorem (SSLT) for dimension one \cite{JinKorepin}
 and also for the higher dimensional case \cite{GK}, 
see Remark~\ref{SSLT} below. More precisely,
the asymptotic behavior of EE
as the size of the subsystem of interest becomes large
for the $XX$ spin chain model with a transverse magnetic field
was analyzed rigorously in \cite{JinKorepin}
using a certain Fisher--Hartwig theorem for large Toeplitz
determinants established in \cite{Basor},
see also \cite{BasorTracy,Ehrhardt}
(self-correlations for a translation invariant spin chain
can be expressed in terms of Toeplitz determinants).
For the more general $XY$ model with a transverse magnetic field,
the authors in \cite{Its} used the Riemann--Hilbert approach and
the steepest descent method to find an explicit expression for EE
(in the asymptotic regime where $XY\to XX$, the results in \cite{JinKorepin}
are recovered from the corresponding 
expression in \cite{Its}).

Note that the above mentioned results are for the one dimensional case only
(we refer the reader e.g.~to the references in \cite{Its,Plenio} for further results
concerning spin chains and also harmonic chain systems).
Much less is known about the higher dimensional case.
It was shown in \cite{Plenio} that for the harmonic lattice model
the entaglement entropy of a cubic region of size $\lambda$
behaves like $\lambda^{d-1}$, $\lambda\to\infty$, 
where $d>1$ is the dimension, 
i.e.~the EE is of the order of the area of the boundary of the cube. 
In the physics community this type of behavior 
is referred to as the {\em area law} for the entropy. 
The area law was initially discovered in the context of the so-called
{\em geometric entropy,\ }which is a component 
of the Bekenstein--Hawking black hole entropy, see \cite{Bombelli,Srednicki}
and the references in \cite{Plenio}.

It is known that there is a correspondence between
$1D$ spin models and systems of non-interacting
fermions on a $1D$ lattice by means of the Wigner--Jordan
transformation. It is therefore natural to consider EE for fermionic
systems in higher dimensions.
It was noticed recently \cite{GK,Wolf},
 that the EE for a system of free fermions
of  arbitrary dimension 
 on a lattice or in the continuum
at zero temperature 
violates the area law and is, in particular,
of larger order than the boundary area of the region in which the
entropy is evaluated. In \cite{GK}
 the authors make a connection
between the asymptotics of EE and a conjecture 
due to Widom \cite{Widomoperator,W5} (see \eqref{Wconj} below), 
and then utilize this connection to posit a formula for EE of a continuous
system, see \eqref{eqconj} below.
The conjecture of Widom appeared originally in the
context of time--frequency limiting problems,
i.e.~problems that involve the extraction of information about a signal
from a measurement in a finite time---finite frequency window.

In \cite{GK}, some of the results are presented without proof.
In this paper, we study Szeg\"o type asymptotics for 
operators of pseudodifferential type with discontinuous symbols
in the higher dimensional case: Various specializations of 
these asymptotics provide the proofs of most, but not all,
of the results left open in \cite{GK}, as explained below.

To fix ideas we recall first some basic facts concerning 
EE (see e.g.~\cite{BBPS}).
Let $H_A,H_B$ be two 
Hilbert spaces, which we assume for simplicity to be finite dimensional.
Using the Schmidt decomposition \cite{Peres},
any state $\psi\in H_A\otimes H_B$
can be expressed in the form $\psi=\sum_i c_i
\phi_{A,i}\otimes \phi_{B,i}$, 
where $0<c_i\leq1$, $\sum_i c_i^2=1$ and $\phi_{A,i},\phi_{B,i}$ 
are orthonormal in $H_A,H_B$,
respectively. 
Associate with $\psi$ the density matrix $\rho=\psi\psi^*$
 acting on $H_A\otimes H_B$. 
 The reduced density matrix $\rho_A$ acting 
on $H_A$ is defined as 
$\rho_A=\Tr_{H_B}\rho=\sum_i c_i^2\phi_{A,i}\phi_{A,i}^*$
and similarly $\rho_B=\Tr_{H_A}\rho=\sum_i c_i^2\phi_{B,i}\phi_{B,i}^*$.
It is easy to check that the matrices $\rho_A$ and $\rho_B$
are well-defined independent of the choice of 
the orthonormal vectors $\phi_{A,i}$, $\phi_{B,j}$.

 The entanglement entropy of a state $\psi\in H_A\otimes H_B$ 
measures how far the state
$\psi$ is from a product state of the form $\phi_A\otimes\phi_B$,
and is defined as 
the von Neumann entropy of {\it either\ }of the reduced density matrices
$$
   S\equiv-\Tr(\rho_A\log_2\rho_A)=-\Tr(\rho_B\log_2\rho_B)
                 =-\sum_i c_i^2\log_2 c_i^2
$$
which is precisely the Shannon entropy of the (squared) 
Schmidt coefficients $c_i$.

In many problems, one is interested in finding EE for $\psi$
which is the ground state of some general many body system.
In \cite{klich} the EE of a system 
of non-interacting (free) fermions in the ground state
was studied. 
We note that the case of a general system with interactions
is very difficult and at present time very little seems to be known
for dimensions higher than one
(see however \cite{Hamma} and \cite{Plenio} where the area law is derived
for the Kitaev model and for the harmonic lattice model, respectively).

Let $\Gamma\subset\R^d$ (resp., $\Gamma\subset\mathbb{T}^d$) 
denote a compact set in momentum space
and fix a compact region $\Omega\subset\R^d$ (resp., $\Omega\subset\Z^d$) 
in position 
space for the continious
(resp., lattice) case.
The ground state of the
system is defined by the projection $P$ 
in $L^2(\R^d)$ (resp., $l^2(\Z^d)$)
onto the modes in the Fermi sea $\Gamma$.
We study the entropy of entanglement between 
fermions located in a 
compact region $\Omega$, scaled by some large $\lambda$,
in position space, and its complement.
Let $Q=\chi_{\lambda\Omega}$ be the projection 
onto $\lambda\Omega$ in $L^2(\R^d)$ (resp., $l^2(\Z^d)$). 

For the system at hand, all the important quantities
for the entanglement problem
can be described in terms of the operator $PQP$ \cite{klich};
in particular the average number of fermions in $\lambda\Omega$ 
is given by ${<}N{>}=\Tr PQP$, the particle number variance
is 
\be\la{fnvar}
   (\Delta N)^2=\Tr [PQP(1-PQP)]
\ee
and the EE is given by 
\begin{eqnarray}\label{entropy}
   S\equiv S_{\Omega,\Gamma}(\lambda)=\Tr h(PQP)
\end{eqnarray}
where 
\be\la{eqh}
   h(t)=-t\log_2 t -(1-t)\log_2(1-t).
\ee
Note that $P=\F\chi_{\Gamma}\F^{-1}$ where $\F$
denotes either the Fourier transform or the Fourier series in the continuous,
lattice case, respectively.

Assuming the applicability of the Widom conjecture
\eqref{Wconj} below to the function \eqref{eqh}
the authors in \cite{GK} suggest 
the  following explicit leading order asymptotics for the EE of
a continuous system as $\lambda\to\infty$
\be\la{eqconj}
     S_{\Omega,\Gamma}(\lambda) 
      = \frac{{\lambda}^{d-1}\log_2{\lambda}}{(2\pi)^{d-1}} \,\frac1{12}
                  \int_{\partial\Omega}\int_{\partial\Gamma} |n_x\cdot n_{p}| \,
                            dS_x\,dS_\xi+ o({\lambda}^{d-1}\log_2{\lambda}),
\ee
where $n_x,n_\xi$ are outward unit 
normals to the (smooth) boundaries $\partial\Omega,\partial\Gamma$
and $dS_x,dS_\xi$ are the area elements.
Although the formula \eqref{eqconj}
was conjectured in \cite{GK}
only for continuous models with smooth
boundaries, it is probably also true for piecewise smooth
boundaries in both the continuous and the lattice
case. Indeed, \eqref{eqconj} was recently checked numerically 
in the lattice case
for $d=2$ in \cite{German}
and for $d=2,3$ in \cite{USC}, and an extremely close
agreement concerning both the order and the (leading)
coefficient was found.

The following two results proven in \cite{GK} 
provide corroborating evidence towards \eqref{eqconj}.
For the cubic domains $\Gamma=[-\frac12,\frac12]^d$
and $\Omega=[0,1]^d$ (resp., $\Omega=\{0,1\}^d$) in the continuous
 (resp., lattice)
case, the following holds as $\lambda\to\infty$
\be\la{eetwosided}
   \frac1{2}\bigg(\frac{\lambda}{2\pi}\bigg)^{d-1}
     S_{1}(\lambda)
         \leq S_{\Omega,\Gamma}(\lambda)
      \leq  d\,\bigg(\frac{\lambda}{2\pi}\bigg)^{d-1}
          S_{1}(\lambda)
\ee
where $S_{1}(\lambda)$ is the entanglement entropy
for the one dimensional system with $\Gamma_1=[-\frac12,\frac12]$
and $\Omega_1=[0,1]$ (resp.,~$\{0,1\}$).
We note that \eqref{eetwosided} together with a result
in \cite{JinKorepin} for the lattice case
\be\la{eqJiKo}
         S_{1}(\lambda) =\frac13\log_2\lambda + 
              o(\log_2\lambda),\qquad\lambda\to\infty,
\ee
implies
\be\la{GKW}
   \frac1{6}\bigg(\frac{\lambda}{2\pi}\bigg)^{d-1}
     \log_2\lambda
         \leq S_{\Omega,\Gamma}(\lambda)
      \leq  \frac{d}{3}\,\bigg(\frac{\lambda}{2\pi}\bigg)^{d-1}
           \log_2\lambda
\ee
thereby demonstrating the violation of the entropy area law 
for {\em cubic domains in the lattice case\ }in a way which is consistent
with \eqref{eqconj}.

In \cite{GK} the authors also prove that for arbitrary 
(measurable) compact $\Omega,\Gamma$ in both the
continuous and lattice cases
\be\la{eqestnv}
        4(\Delta N)^2 \leq 
           S_{\Omega,\Gamma}(\lambda) 
               \leq C\cdot(\log_2\lambda)\cdot(\Delta N)^2
\ee
where the constant $C$ depends only on the dimension $d$.
(The estimate \eqref{eqestnv} for the {\em lattice\ }case with $d=1$
was proved in \cite{Fannes}. The proof of \eqref{eqestnv}
in \cite{GK}
for all $d\geq1$
in the lattice case is analogous, but the continuous case requires a new idea
as provided in \cite{GK}.)

The fact that Theorem \ref{CH3_lem.lap.1} below
is sharp implies, together with \eqref{fnvar}, 
that in the continuous case
if $\Gamma,\Omega$ have $C^1$ boundaries then for some $c_1,c_2>0$
\be\la{thisp1}
     c_1 \lambda^{d-1}\log_2\lambda \leq (\Delta N)^2 \leq
                c_2 \lambda^{d-1}\log_2\lambda,\qquad \lambda\to \infty,
\ee
and also that for any $\beta\in(0,1)$ there exists (a Cantor-like) set $\Gamma$
such that for $\Omega=[0,1]^d$ for some $c_1,c_2>0$
\be\la{thisp2}
     c_1 \lambda^{d-\beta} \leq (\Delta N)^2 \leq
                c_2 \lambda^{d-\beta},\qquad \lambda\to \infty.
\ee
The inequalities \eqref{eqestnv}, \eqref{thisp1} yield
the following:
In the continuous case,
if $\Omega,\Gamma$ have $C^1$ boundaries then for
some $c_1,c_2>0$ that depend on $\Omega,\Gamma$
\be\la{eq1estnv}
        c_1\lambda^{d-1}\log_2\lambda 
      \leq S_{\Omega,\Gamma}(\lambda) \leq 
          c_2\lambda^{d-1}(\log_2\lambda)^2,\qquad
                                    \lambda\to\infty.
\ee
Note that \eqref{eq1estnv} proves the violation of the area law
in the {\em continuous case for arbitrary domains\  }with $C^1$ boundary
and also gives the expected order for the lower estimate 
(and also the expected upper estimate, up to a power of $\log_2\lambda$)
consistent with \eqref{eqconj}.
\begin{remark}
It was proved independently in \cite{Wolf} that for 
cubic domains in the lattice case the estimate \eqref{eq1estnv} holds.
Note that in this case we have the stronger estimate \eqref{GKW}.
\end{remark}
Concerning irregular boundaries, \eqref{eqestnv} 
together with \eqref{thisp2}
shows that in the {\em continuous\ }case
in any dimension $d$ and
for any $\beta\in(0,1)$, there exists a 
set $\Gamma$ so that for $\Omega=[0,1]^d$
for some $c_1,c_2>0$
\be\la{estimatefractal}
           c_1 \lambda^{d-\beta} \leq
                 S_{\Omega,\Gamma}(\lambda)\leq
                  c_2 \lambda^{d-\beta}\log_2\lambda,
     \qquad \lambda\to\infty.
\ee
Finally, we note that in the one dimensional 
lattice case the following result is proved in \cite{Fannes}:
For any $\beta_1\in(0,1)$ there 
exists a set $\Gamma_1\subset{\mathbb T}$ 
such that for $\Omega_1=\{0,1\}$
for some $c_1,c_2>0$
$$
         c_1 \lambda^{1-\beta_1} \leq 
                S_{\Omega_1,\Gamma_1}(\lambda) \leq
                c_2 \lambda^{1-\beta_1}\log_2\lambda, \qquad \lambda\to \infty.
$$
Combining this result with 
a straightforward generalization of \eqref{eetwosided}
to the case $\Gamma=\Gamma_1^d$, $\Omega=\Omega_1^d$,
we see that \eqref{estimatefractal} holds also for the {\em lattice\ }case
(where $\beta=1-\beta_1$)\footnote{The fractal set $\Gamma$ in \cite{Fannes}
is very similar to the fractal set appearing in Lemma~\ref{lemex} below,
which was constructed in earlier work of the author,
see math.FA/0212215, math.CA/0212254.}.

The outline of the paper is as follows: In Section~2, we describe
the main results, the proofs of these results together with
various technical lemmas are given in Section~3.
\section{Main results}
We now state the main results of the paper.
As noted above, some of these results
were used in \cite{GK} to estimate EE
for certain physical systems.

The problem of evaluating the Szeg\"o type asymptotics 
for operators of pseudodifferential type with
symbols discontinuous in both position and momentum
 in the higher dimensional case
was introduced in \cite{W5} (see Remark~\ref{remconj}
below),
we refer also to \cite{Wevd,LS} where certain related results
can be found. See e.g.~\cite{W3} for a more detailed account
of the results in this area.

Let $d\in\N$, let $\meas$ denote the Lebesgue measure
in $\RD$, set $H=L^2(\RD)$ 
and let $\|\cdot\|_{k}$ be the standard norm in $L^k(\RD)$,
$k=1,2,\infty$.
Denote by $\|\cdot\|$, $\|\cdot\|_{\HS}$, $\|\cdot\|_{\mathfrak{S}_1}$
the operator norm, the Hilbert--Schmidt and the trace-class norm in $H$,
respectively.
Let ${\mathcal{S}}(\RD)$ stand for the Schwartz space
and denote by $\int$ the integration over $\RD$.
Let $\F$ and $\,\hat{}\,$ denote the Fourier transform:
$
        \F_{x\ra{}u}[{g}(x)] \equiv \hat{g}(u) 
   =\int e^{-iu\cdot{}x}{g}(x)\,dx
$, ${g}\in{}H$.
For a function ${g}\in{}L^2(\RD)$ denote its $L^2$ modulus
of continuity by
$
    \omega_2[{g}](h)= \|{g}(\cdot+h)-g(\cdot)\|_2
$, $h\in\RD$.
For a set $\Omega\subset\RD$ denote by $\chi_\Omega$
its characteristic function.
We characterize the regularity of the set $\Omega$
in terms of $\omega_2[\chi_\Omega]$:
assume that $\chi_\Omega\in{}L^2(\RD)$
(i.e.~$\meas(\Omega)<\infty$)
and that there exist
$0<\beta_\Omega\leq1$ and $c_\Omega>0$ 
such that for small enough $|h|$
\begin{equation}
\label{eq_charos}
        (\omega_2[\chi_\Omega](h))^2\leq c_\Omega|h|^{\beta_\Omega}
\end{equation}
where $|h|=(\sum_{j=1}^d|h_j|^2)^{1/2}$.
Let $\Omega_h=\{x-h|x\in\Omega\}$, $h\in\RD$. 
The left-hand side of \eqref{eq_charos} equals 
\begin{equation}
\label{eq_note}
    \|\chi_\Omega(\cdot+h)-\chi_\Omega(\cdot)\|_{2}^2 
  = \meas(\Omega\setminus\Omega_{h}) 
     + \meas(\Omega\setminus\Omega_{-h}),
\end{equation}
which gives the geometrical meaning of \eqref{eq_charos}.
Introduce a projection $P:H\to{}H$ by
$(P{g})(x)=\chi_\Omega(x){g}(x)$, $x\in\RD$, ${g}\in{}H$.
Following \cite{W5} we 
consider a family 
of integral operators $A_\lambda:H\to H$, $\lambda\geq2$, 
 of pseudodifferential type 
with a non-smooth or discontinuous symbol.
Let $A_\lambda$, $\lambda\geq2$, 
have the kernel 
\begin{equation}
\label{eqsaf4}
   K_{{A_{\lambda}}}(x,y) = \bigg(\frac{\lambda}{2\pi}\bigg)^d
              \int e^{i\lambda\xi\cdot{}(x-y)} 
                                  \sigma(x,y,\xi)\,d\xi 
\end{equation}
where $\sigma$ (is measurable and) satisfies the following mild condition:
Define 
\begin{equation}
\label{eqphi}
     \phi(u)=\sup_{x,y}\big|\F_{\xi\to u}[\sigma(x,y,\xi)]\big|,
               \qquad \psi(u)=\phi(u)\phi(-u)
\end{equation}
and assume $\psi\in L^{1}(\RD)$ and
for certain $0<\beta\leq1$ and $c>0$ 
\begin{equation}
\label{eqsaf1}
     \int_{|u|\geq\rho}\psi(u)\,du\leq{}c\rho^{-\beta},\qquad\rho\geq1.
\end{equation}
See Remark \ref{remex} below for an example of such $\sigma$.
Setting $\lambda=h^{-1}$, $h>0$, in \eqref{eqsaf4} we can obtain
semiclassical type asymptotics.
 We study the asymptotics of the trace of 
\begin{equation}
\label{eqfirst}
      Pf(PA_\lambda{}P)P- Pf(A_\lambda{})P,
\end{equation} 
as $\lambda\ra\infty$, for suitable functions $f$. 
For the particular 
choice $\sigma(x,y,\xi)=\chi_\Gamma(\xi)$, $\Gamma\subset\RD$,
such a study is motivated by the following question:
what can be said about a function if its restriction to 
$\Omega$ and the restriction 
of its Fourier transform to $\Gamma$ are known, see \cite{W5}.
An order sharp estimate for the trace of
\eqref{eqfirst} 
is found 
in Theorem~\ref{CH3_s0cor1} below in two settings; 
we assume either that $f$ is analytic,
or $A_\lambda$ is self-adjoint and $f$ has a bounded second derivative.
Theorem~\ref{CH3_s0cor1}
is a generalization of the classical Szeg\"o limit theorem \cite{Sz1}, see
e.g., \cite{W3,LS} for a review of related results.
Let $\log$ denote the natural logarithm.
The following result is basic for the proof of Theorem~\ref{CH3_s0cor1}.
\begin{theorem}
\label{CH3_lem.lap.1}
Let $\beta$, $\sigma$, $A_\lambda$, $\beta_\Omega$, 
$\Omega$ be as above.
Then
there exist two constants $c(\Omega,\sigma)$ and $C(\Omega,\sigma)$ 
such that one has, for $\lambda\geq2$,
\begin{equation}
\label{CH3_s1e2222}
\begin{aligned}
        \|P{}{}{A_{\lambda}}\|_{\HS}^2
            \leq c(\Omega,\sigma)\cdot\lambda^{d}, 
\end{aligned}
\end{equation}
and
\begin{equation}
\label{CH3_s1e2223}
\begin{aligned}
        \|P{}{}{A_{\lambda}}(I-P{})\|_{\HS}^2
            \leq C(\Omega,\sigma)\cdot
                     \begin{cases}
            \lambda^{d-\min(\beta,\beta_\Omega)},&\beta\neq\beta_\Omega\cr
                                 \lambda^{d-\beta}\log\lambda,&\beta=\beta_\Omega.
                     \end{cases}
\end{aligned}
\end{equation}
and the same estimate holds for $\|(I-P){A_{\lambda}}P\|_{\HS}^2$.
The estimate \eqref{CH3_s1e2223} is {\em sharp\ }on 
the described class of $\sigma$, 
that is for a certain $\sigma$ (which can be chosen so that 
the corresponding $A_\lambda$ is self-adjoint),
the reverse inequality to \eqref{CH3_s1e2223} holds 
for some (different) constant $C(\Omega,\sigma)$.
\end{theorem}

For an analytic on some disc $f(z)=\sum_{m=0}^\infty{}c_m{}z^m$
set
\begin{equation}
\label{CH3_s0eqfstar}
          {}f_*(z) = \sum_{m=2}^\infty{}m(m-1)|c_m|{}z^{m-2}
\end{equation}
and
\begin{equation}
\label{CH3_s0esup}
               S(\sigma)= \sup_{\lambda\geq2}\|A_{\lambda}\|.
\end{equation}
The condition $S(\sigma)<\infty$ holds for a wide class of $\sigma$,
see Remark \ref{remex} below.
\begin{theorem}
\label{CH3_s0cor1}
Let the assumptions of Theorem~\ref{CH3_lem.lap.1} hold.
Assume further that either 

(i) $S(\sigma)<\infty$ and $f$ is analytic
on a neighborhood of 
$\{z:|z|\leq S(\sigma)\},$
or 

(ii) $A_{\lambda}$ is self-adjoint for all $\lambda\geq2$, 
and  $f$ is such that
$f^{\prime\prime}\in{}L^\infty({\mathbb{R}})$. 

\noindent
 Then the operator \eqref{eqfirst}
is trace-class. Moreover, there exists a constant $C(\Omega,\sigma)$ such that
one has, for $\lambda\geq2$,
\begin{equation}
\label{CH3_s0e17four}
    \Big|\tr\big[Pf(P{}A_{\lambda}P)P -
           Pf(A_{\lambda})P\big]\Big|
   \leq{}C(\Omega,\sigma)\cdot{}\tilde{C}(f)
    \cdot\begin{cases}
                  \lambda^{d-\min(\beta,\beta_\Omega)},&\beta\neq\beta_\Omega\cr
                                 \lambda^{d-\beta}\log\lambda,&\beta=\beta_\Omega,
                  \end{cases}
\end{equation}
where
$$
      \tilde{C}(f)=\frac12\cdot
                 \begin{cases}f_*\big(S(\sigma)\big),&\text{in case (i)}\cr
                 \|f^{\prime\prime}\|_{\infty},&\text{in case (ii)}.
         \end{cases}
$$
In both cases, the estimate \eqref{CH3_s0e17four} is {\em sharp,\ }that 
is for certain $f$
and $\sigma$, the reverse inequality 
to \eqref{CH3_s0e17four} holds for some (different)
constant $C(\Omega,\sigma)$.
\end{theorem}
In the case of analytic $f$
we prove slightly more, namely that the sharp
estimate \eqref{CH3_s0e17four} holds 
with the trace-class norm of \eqref{eqfirst} in the left-hand side.
See subsections \ref{subsec.2.1} and \ref{subsec.2.2}
 for the proofs of Theorem \ref{CH3_lem.lap.1}
and Theorem \ref{CH3_s0cor1}, respectively.
It is possible to compute 
the leading term in the asymptotics
of $\tr{}Pf(PA_\lambda{}P)P$ in certain special cases.
\begin{corollary}
\label{CH3_s0cor5}
(i) Let $\Omega$ and $\Gamma$ be two bounded
domains with $C^\infty$ boundaries.
Assume that the symbol $\sigma$ does not depend on $y$,
and that $\sigma(x,\xi)=\tau(x,\xi)\chi_\Gamma(\xi)$, 
where $\tau\in{\mathcal{S}}(\R^{2d})$.
Let $f$ be analytic on a neighborhood of the disc
$
              \{z:|z|\leq\|\tau\|_{\infty}\}
$ 
and satisfy $f(0)=0$. Then the operator $f(PA_\lambda{}P)$ is trace-class,
and furthermore, for any small enough
$\eps>0$ there exist two constants $C(\Omega,\Gamma,\tau,\eps)$ 
and $\Lambda(\eps)\geq2$ so that, for $\lambda\geq\Lambda(\eps)$,
\begin{equation}
\label{CH3_s0e17prim}
\begin{aligned}
    \bigg|\tr f(PA_\lambda{}P) &- 
         \bigg(\frac{\lambda}{2\pi}\bigg)^d\int_\Gamma\int_\Omega 
              f(\tau(x,\xi))\,dx\,d\xi\bigg|\\
         &\leq{}C(\Omega,\Gamma,\tau,\eps)\cdot{}f_*\big((1+\eps)\|\tau\|_\infty\big)
                     \cdot\lambda^{d-1}\log\lambda.
\end{aligned}
\end{equation}

(ii)
Let \eqref{eq_charos} hold and 
assume that the symbol $\sigma$ is real-valued, 
depends only on $\xi$ and \eqref{eqsaf1} holds.
Let $f^{\prime\prime}\in L^\infty([-\|\sigma\|_\infty,\|\sigma\|_\infty])$ and
assume $f(\sigma(\xi))\in{}L^1(\RD)$.
Then $f(PA_\lambda{}P)$ is trace-class
and there exists a constant $C(\Omega,\sigma)$ 
such that, for $\lambda\geq2$,
\begin{equation}
\label{CH3_s0e17bis}
\begin{aligned}
    \bigg|\tr (Pf(PA_\lambda{}P)P) &- \bigg(\frac{\lambda}{2\pi}\bigg)^d
         \,{\textrm{\em mes}}(\Omega)\,\int
               f(\sigma(\xi))\,d\xi\bigg|\\
         &\leq{}C(\Omega,\sigma)\cdot\|f^{\prime\prime}\|_{\infty}
    \cdot \begin{cases}\lambda^{d-\min(\beta,\beta_\Omega)},
                                                      &\beta\neq\beta_\Omega\cr
                                 \lambda^{d-\beta}\log\lambda, &\beta=\beta_\Omega.
                  \end{cases}
\end{aligned}
\end{equation}
\end{corollary}
The proof of part (i) follows from 
the functional calculus
results developed in \cite{W5}, see subsection \ref{subsec.2.3}. 
Part (ii) is proved as follows. 
Under the above assumptions $A_{\lambda}$ is self-adjoint, $\lambda\geq2$,
 and also the operator $f(A_\lambda)$ is well-defined and has the kernel
$
   K_{f(A_{\lambda})}(x,y) = \big(\frac{\lambda}{2\pi}\big)^d
              \int e^{i\lambda\xi\cdot{}(x-y)} \,
                                  f(\sigma(\xi))\,d\xi.
$
The operator $Pf(A_\lambda)P$, $\lambda\geq2$, is trace-class
since it is a composition of two Hilbert--Schmidt operators
$P\F^{-1}|f|^{1/2}\F$, $\F^{-1}|f|^{1/2}(\sgn{f})\F P$. 
Note also that $K_{f(A_{\lambda})}(x,y)$ is continuous, 
since $f(\sigma(\xi))\in{}L^1(\RD)$.
Now we simply write the trace of $Pf(A_{\lambda})P$
as the integral of its kernel over the diagonal.
\begin{remark} 
Let $\chi_\Omega,\chi_\Gamma$ satisfy \eqref{eq_ts} below
with $0<\beta_\Omega,\beta\leq1$, respectively,
and set $\sigma(x,y,\xi)=\chi_\Gamma(\xi)$ (this example is considered in
the proof of Theorem \ref{CH3_lem.lap.1}).
The order sharp remainder estimate in  \eqref{CH3_s0e17bis}
for $\beta_\Omega\neq\beta$
shows that the set with less regular boundary
contributes to the order of the remainder.
In the case $\beta_\Omega=\beta$
the logarithmic factor
persists even for $\Omega,\Gamma$ with $C^\infty$ boundaries
(in which case $\beta_\Omega=\beta=1$)
due to the fact that the symbol
$
        \chi_\Omega(x) \,\chi_{\Gamma}(\xi)\, \chi_\Omega(y)
$
has discontinuities in {\em both\ }the position variables $x,y$
and the momentum variable $\xi$. (This should be compared with
the power type asymptotics in e.g.~\cite{LS,GO,GiCPDE} when a discontinuity
in momentum only is present.)
\end{remark}
\begin{remark}
The leading term 
$
       \tr Pf(A_{\lambda})P
$
in the asymptotics of 
$\tr{}Pf(PA_{\lambda}P)P$ in \eqref{CH3_s0e17prim} and 
\eqref{CH3_s0e17bis} 
is of Weyl type.
It is written as an integral over the diagonal $(x,x,\xi)$
of the corresponding to $\lambda$ phase volume.
\end{remark}
\begin{remark}
\label{remconj}
In \cite{Widomoperator,W5} the following
 {\em second\ }order generalization 
of \eqref{CH3_s0e17prim} was conjectured:
Under the assumptions of Corollary \ref{CH3_s0cor5}(i)
 the following holds, as $\lambda\to\infty$
\begin{equation}
\label{Wconj}
\begin{aligned}
 \Tr f(PA_\lambda P) &= \Big(\frac{\lambda}{2\pi}\Big)^d\,
       \int_\Omega\int_\Gamma f(\tau(x,\xi))\,dxd\xi
       \\  &+ \Big(\frac{\lambda}{2\pi}\Big)^{d-1}
               \frac{\log{\lambda}}{4\pi^2}
     \int_{\partial\Omega}\int_{\partial\Gamma}
          |n_x\cdot n_{p}| \,U(0,\tau(x,\xi);f) \,dS_x d S_{\xi}\\ &+ o({\lambda}^{d-1}\log{\lambda})
\end{aligned}
\end{equation}
where $n_x,n_{\xi}$ are the outward unit normals 
to $\partial\Omega,\partial\Gamma$, respectively, and 
$$
   U(a,b;f) =\int_0^1 \frac{f((1-t)a+tb)-[(1-t)tf(a)+tf(b)]}{t(1-t)}\,d t.
$$ 
This conjecture is still open and it was one of the motivations
for the present work.
\end{remark}
\begin{remark}
\label{SSLT}
In a broader context, the Widom conjecture \eqref{Wconj} is
 a generalization of the
strong (two-term) Szeg\"o limit theorem (SSLT) for 
the {\em continuous\ }setting. (As noted earlier, a generalization
of the SSLT plays a
central role in the computation of the EE for the $XX$ spin
chain model in \cite{JinKorepin}.)
The SSLT was initially used by Onsager in his celebrated computation of the spontaneous magnetization for
the $2D$ Ising model (see e.g.~\cite{Bottcher} and the references therein).
 It is interesting to note that in
 Onsager's computation (and also in \cite{JinKorepin}) the leading
asymptotic term vanishes, and one needs to compute the sub-leading term.
This is exactly the
situation in \cite{GK}:
the leading
term should vanish since for the entropy 
function $h$ in \eqref{eqh}, $h(1)=0$ (and also $h(0)=0$ holds).
\end{remark}
\begin{remark}
\label{remex}
Let $\Gamma\subset\RD$ be
such that for some $0<\beta\leq1$ the function
$\chi_\Gamma$ satisfies \eqref{eq_charos}
with $\gamma=\beta$ (see also Lemmas \ref{lemex}, \ref{CH3_s3.lem.1} below).
Set $\sigma(x,y,\xi)=\tau(x,y,\xi)\,\chi_\Gamma(\xi)$, where 
$\tau$ is satisfies for some $c=c(\tau)<\infty$
\begin{equation}
\label{eqsaf50}
          \sup_{x,y}\big| \F_{\xi\ra{}u}[\tau(x,y,\xi)]\big|\leq
c\,(1+|u|)^{-d-1},\qquad u\in\RD.
\end{equation}
A standard application of 
the Cauchy inequality implies that there is $C=C(\tau,d)<\infty$ so that
$
    \psi(u)\leq C \int|\hat\chi_\Gamma(u-v)|^2(1+|v|)^{-d-1}dv\in L^1(\RD),
$
and \eqref{eqsaf1} holds
(see, for instance, \cite[Section~3.4.2]{Gi}).
For $\sigma(x,y,\xi)=\chi_{\Gamma}(\xi)$ we can take 
$\psi(u)=|\hat\chi_\Gamma(u)|^2$.
Now for \eqref{CH3_s0esup}. 
If $\sigma(x,y,\xi)$ is a classical zeroth 
order (parameter dependent) 
symbol in the sense of pseudodifferential operators,
then $S(\sigma)<\infty$ is a standard result, 
see e.g.~\cite{GS,W5}. 
However we are interested here
in $\sigma$ with limited regularity for which no such general
results are available. We restrict ourselves to 
the following standard example in which
irregularities in $x$ and, most importantly, in $\xi$
are allowed.
Assume $\sigma=\sigma(x,\xi)$ (we only need $S(\sigma)<\infty$
when $A_\lambda$ is not assumed to be self-adjoint)
is of the form $\tau(x,\xi)\,\chi(\xi)$.
Then $A_\lambda=\tilde{A}_\lambda\F^{-1}\chi\F$
where $\tilde{A}_\lambda$ has integral kernel
$\big(\frac{\lambda}{2\pi}\big)^d\int e^{i\lambda(x-y)\cdot\xi}\tau(x,\xi)\,d\xi$.
Assume that $\chi\in L^\infty(\RD)$, then $\|\F^{-1}\chi\F\|\leq\|\chi\|_\infty$
(the factor $\chi$ is allowed to be quite irregular).
Now $\sup_{\lambda\geq2}\|\tilde{A}_\lambda\|<\infty$
follows from the standard estimate for the bilinear form under 
the sole assumption
 that $\tilde{\phi}(u)=\sup_x|\F_{\xi\to u}\tau(x,\xi)|\in L^1(\RD)$.
Indeed, by the Cauchy inequality for any $f,g\in{}H$
$$
\begin{aligned}
    (&2\pi)^d \big|(\tilde{A}_\lambda f,g)\big|
              \leq \lambda^d\int\!\!\!\int \tilde{\phi}(-\lambda(x-y))
                         \,|f(y)|\,|g(x)|\,dx\,dy\\
          &\leq \Big(\lambda^{d}\!
                           \int\!\!\!\int \tilde{\phi}(-\lambda(x-y))
              \,|f(y)|^2dxdy\Big)^{1/2}
 \Big(\lambda^{d}\!
                           \int\!\!\!\int \tilde{\phi}(-\lambda(x-y))
\,|g(x)|^2dxdy\Big)^{1/2}\\
            &\leq \|\tilde{\phi}\|_1\, \|f\|_2\,\|g\|_2
\end{aligned}
$$
which implies 
$      S(\sigma)
            \leq  (2\pi)^{-d}\,\|\tilde{\phi}\|_1$. 
Note that if $\sigma(x,\xi)=\tau(x,\xi)\chi_\Gamma(\xi)$
where $\Gamma$ is as before,  and
$\sup_x|\F_{\xi\to u}\tau(x,\xi)|\leq{}c\,(1+|u|)^{-d-1}$
as in \eqref{eqsaf50}, then both \eqref{eqsaf1} 
and $S(\sigma)<\infty$ hold.
\end{remark}
In the proof of Theorem \ref{CH3_lem.lap.1} 
we use the following auxiliary results.
They are concerned with a {\em two-sided} 
version of the estimate \eqref{eq_charos}: assume that 
  ${g}\in{}L^2(\RD)$ 
is such that there exist $0<\gamma<2$
and $c_1,c_2>0$ so that for small enough $|h|$
\begin{equation}
\label{eq_ts}
       c_1 |h|^{\gamma}\leq (\omega_2[{g}](h))^2 \leq c_2|h|^{\gamma}.
\end{equation}
\begin{lemma}
\label{lemex}
  For any $d\in\N$ and $0<\beta_\Omega\leq1$ there exists
a compact 
set $\Omega\subset\RD$ such that 
$\chi_\Omega$ satisfies
\eqref{eq_ts} with $\gamma=\beta_\Omega$.
\end{lemma}
\begin{lemma}
\label{CH3_s3.lem.1} 
Assume that a function $f\in{}L^2(\RD)$
satisfies \eqref{eq_ts} for some $0<\gamma<2$ and $c_1,c_2>0$. 
Then

(i) there exist $b_1,b_2>0$ such that
\begin{equation}
\label{eqnewww}
          b_1\,\rho^{-\gamma} \leq 
                   \int_{|\xi|\geq{}\rho} 
                      |\hat{f}(\xi)|^2\,d\xi 
           \leq b_2\,\rho^{-\gamma},\qquad\rho\geq1;
\end{equation}

(ii) the upper estimate in \eqref{eq_ts}
implies the upper estimate in \eqref{eqnewww}.
\end{lemma}
Lemma \ref{lemex} is proved in subsection \ref{subsec.2.4}.
Lemma  \ref{CH3_s3.lem.1} for $\gamma=1$
was proved in \cite[Lemma~2.10, 4.2]{BCT},
the proof for $\gamma\in(0,2)$ is analogous and is left to the reader,
see also \cite{Gi2} and \cite[Lemma~3.4.1]{Gi}.
(If one introduces an average of $\omega_2[{g}](h)$
over $\|h\|_d\leq\eps$, 
then the upper estimates in this modification of \eqref{eq_ts},
and in \eqref{eqnewww}
become {\em equivalent,\ }and so do the lower ones, see \cite{Cl,Gi2}.)
\section{Proofs}
\subsection{Proof of Theorem \ref{CH3_lem.lap.1}}
\label{subsec.2.1}
We denote $\lambda$-independent constants 
(that may depend on $\Omega,\sigma$)
by $c_k$, $k\in\N$.
Let us consider $\|P{}{}{A_{\lambda}}(I-P{})\|_{\HS}$ only,
the case of $\|(I-P){A_{\lambda}}P)\|_{\HS}$ is completely analogous.

1. Let us prove first the upper estimate \eqref{CH3_s1e2223}.
Using \eqref{eqphi} we obtain
\begin{equation}
\label{eqsaf10}
\begin{aligned}
      \null&\|P{}{}{A_{\lambda}}(I-P{})\|_{\HS}^2
                                =\lp\frac{\lambda}{2\pi}\rp^{2d}
          \int\!\!\int\!\!\int\!\!\int e^{i\lambda\xi\cdot(x-y)}\,
                        e^{-i\lambda\eta\cdot(x-y)}\,
                \sigma(x,y,\xi)\,\overline{\sigma(x,y,\eta)} \\
             &\qquad\qquad\qquad\qquad\qquad\qquad\qquad\qquad
                 \qquad{}{}{}\times 
           \,\chi_\Omega(x)\, (1-\chi_\Omega(y))d\xi\,dy\,d\eta\,dx\\
     &=\lp\frac{\lambda}{2\pi}\rp^{2d}
          \int\!\!\int  
  \F_{\genfrac{}{}{0pt}{}{\xi\ra{}-\lambda(x-y)}
                   {\eta\ra{}\lambda(x-y)}}[\sigma(x,y,\xi)\,
                                       \overline{\sigma(x,y,\eta)}]\,\,
         \chi_\Omega(x)\, (1-\chi_\Omega(y))\,dx\,dy\\
   &\leq   \lp\frac{\lambda}{2\pi}\rp^{2d}
          \int\!\!\int \psi\big(\lambda(x-y)\big) \,\chi_\Omega(x)\, [1-\chi_\Omega(y)]\,dx\,dy.
\end{aligned}
\end{equation}
Changing variables $x-y=2x^\prime$, $x+y=2y^\prime$
and dropping the primes we rewrite the right-hand side of \eqref{eqsaf10}
in the form
\begin{equation}
\label{eqsaf19}
\begin{aligned}
       \lp\frac{\lambda}{2\pi}\rp^{2d}
          \int\psi\big(2\lambda{}x\big) \,dx\,
            \int \chi_\Omega(x+y)\, 
                 [1-\chi_\Omega(-x+y)]\,dy. 
\end{aligned}
\end{equation}
Since $\chi_\Omega$ takes on values $0$ or $1$ only, we have
$$
\begin{aligned}
       \int \chi_\Omega(x+y)[1&-\chi_\Omega(-x+y)]\,dy
         =\|\chi_\Omega\|_{2}^2 - \int \chi_\Omega(x+y)\,\chi_\Omega(-x+y)\,dy\\
         &=\frac12\int [\chi_\Omega(x+y)-\chi_\Omega(-x+y)]^2\,dy
         =\frac12\,(\omega_2[\chi_\Omega](2x))^2.
\end{aligned}
$$
Therefore \eqref{eqsaf10} and \eqref{eqsaf19} imply
\begin{equation}
\label{eq_integral}
\begin{aligned}
      \|P{}{}{A_{\lambda}}&(I-P{})\|_{\HS}^2
    \leq \lambda^{2d}\,c_1 \int\psi\big(2\lambda{}x\big) \,(\omega_2[\chi_\Omega](2x))^2\,dx \\
  &=
  \lambda^{2d}\,c_1\bigg(\int_{|x|\leq1/\lambda}+\int_{1/\lambda\leq|x|\leq1}
    +\int_{|x|\geq1}\bigg)
        \, \psi\big(2\lambda{}x\big) \,(\omega_2[\chi_\Omega](2x))^2\,dx.
\end{aligned}
\end{equation}
Using \eqref{eq_charos} and 
making a change of variables $u=2\lambda{}x$ we estimate 
the first integral in \eqref{eq_integral} as follows
\begin{equation}
\label{firsteq}
\begin{aligned}
   \lambda^{2d}\,c_1&\int_{|x|\leq1/\lambda}
        \psi\big(2\lambda{}x\big) \,(\omega_2[\chi_\Omega](2x))^2\,dx
   \leq \lambda^{2d}\,c_2\int_{|x|\leq1/\lambda}
        \psi\big(2\lambda{}x\big) \,|x|^{\beta_\Omega}\,dx\\
    &\leq \lambda^{d-\beta_\Omega}\,c_3
           \int_{|u|\leq2}\psi(u)|u|^{\beta_\Omega}\,du
            \leq \lambda^{d-\beta_\Omega}\,
                        c_4\,\|\psi\|_{1},\quad\lambda\geq2.
\end{aligned}
\end{equation}
Next, noting that 
$
    (\omega_2[\chi_\Omega](2x))^2  \leq
          2\|\chi_\Omega\|_{2}^2 = 2\,\meas(\Omega),
$
setting $u=2\lambda{}x$ 
and using \eqref{eqsaf1} we estimate
 the third integral in \eqref{eq_integral} in the following way
\begin{equation}
\label{thirdeq}
\begin{aligned}
   \lambda^{2d}\,c_1\int_{|x|\geq1}
        &\psi\big(2\lambda{}x\big) \,(\omega_2[\chi_\Omega](2x))^2\,dx
   \leq 2\,\meas(\Omega)\,\lambda^{2d}\,\int_{|x|\geq1}
        \psi\big(2\lambda{}x\big) \,dx\\
    &\leq \lambda^d\,c_5\int_{|u|\geq2\lambda}\psi(u)\,du
    \leq \lambda^{d-\beta}\,c_6,
                                   \quad\lambda\geq2.
\end{aligned}
\end{equation}
Consider now the second integral in \eqref{eq_integral}.
Set
$
  \Psi(r) =\int_{\SD} \psi(r\theta)\,dS_\theta$, $r\geq1,
$
and note that by \eqref{eqsaf1}
\begin{equation}
\label{eqsaf400}
    \int_{r}^\infty \Psi(s) \,s^{d-1}\,ds \leq c\, r^{-\beta},\qquad r\geq1.
\end{equation}
Then the second integral in \eqref{eq_integral} is estimated as follows
\begin{equation}
\label{eqsaf40}
\begin{aligned}
    \lambda^{2d}\,c_1\int_{1/\lambda\leq|x|\leq1}
        \psi\big(2\lambda{}x\big) \,(\omega_2[\chi_\Omega](2x))^2\,dx
          &\leq\lambda^{d-\beta_\Omega}\,c_7
\int_{2\leq|u|\leq2\lambda}\psi(u)\,|u|^{\beta_\Omega}\,du\\
    &= \lambda^{d-\beta_\Omega}\,c_7
           \int_{2}^{2\lambda} \Psi(r)\,r^{\beta_\Omega+d-1}\,dr.
\end{aligned}
\end{equation}
Writing $\Psi(r)=-\frac{1}{r^{d-1}}\frac{d}{dr}\int_{r}^\infty \Psi(s) \,s^{d-1}\,ds$ 
and integrating by parts we obtain
$$
\begin{aligned}
     \lambda^{d-\beta_\Omega}\,c_7 &\int_{2}^{2\lambda}
 \Psi(r)\,r^{\beta_\Omega+d-1}\,dr
   =\lambda^{d-\beta_\Omega}\,c_7\cdot\bigg[ 2^{\beta_\Omega}
 \int_{2}^\infty\Psi(s) \,s^{d-1}\,ds\\
         &-2^{\beta_\Omega}\lambda^{\beta_\Omega}
                   \int_{2\lambda}^\infty\Psi(s) \,s^{d-1}\,ds
        + \beta_\Omega\int_2^{2\lambda}\Big(\int_{r}^\infty\Psi(s) \,s^{d-1}\,ds\Big)
  \,r^{\beta_\Omega-1}\,dr\bigg].
\end{aligned}
$$
Discarding the negative term above and
 using \eqref{eqsaf400} we obtain
$$
\begin{aligned}
     \lambda^{d-\beta_\Omega}\, c_7\int_{2}^{2\lambda} 
               \Psi(r)\,r^{\beta_\Omega+d-1}\,dr
&\leq \lambda^{d-\beta_\Omega}\,c_8
\cdot\Bigg( 1
+ \begin{cases} 
             1,   &\beta_\Omega<\beta \cr
             \log\lambda,                                &\beta_\Omega=\beta \cr
             \lambda^{\beta_\Omega-\beta},     &\beta_\Omega>\beta
    \end{cases}\,\,\Bigg)\\
   &\leq  c_9\cdot \begin{cases} 
             \lambda^{d-\min(\beta,\beta_\Omega)},   &\beta\neq\beta_\Omega \cr
             \lambda^{d-\beta}\,\log\lambda,                &\beta=\beta_\Omega.
    \end{cases}
\end{aligned}
$$
Substituting the latter estimate in \eqref{eqsaf40} 
and collecting the estimates
\eqref{firsteq}, \eqref{thirdeq} we complete the proof of \eqref{CH3_s1e2223}.

2. Let us now prove the sharpness of the estimate \eqref{CH3_s1e2223}.
We choose arbitrary $0<\beta_\Omega,\beta\leq1$.
By Lemma \ref{lemex}
there exist two sets $\Omega,\Gamma$ such that $\chi_\Omega,
\chi_\Gamma\in{}L^2(\RD)$ 
satisfy the two-sided estimate \eqref{eq_ts} with the power
$\beta_\Omega$ and $\beta$, respectively.
Set $\sigma(x,y,\xi)=\chi_\Gamma(\xi)$ (then the corresponding
$A_\lambda$ is self-adjoint) and let $\psi(u)=|\hat{\chi}_\Gamma(u)|^2$.
Apply Lemma \ref{CH3_s3.lem.1}(i)
to $\hat{\chi}_\Gamma$. 
Then for some $b_1,b_2>0$, 
$
      b_1\rho^{-\beta}
          \leq \int_{|u|\geq{}\rho} \psi(u)\,du 
           \leq b_2\rho^{-\beta}
$
for $\rho\geq1$.
Let $\kappa$ be such that $b_2\kappa^{-\beta}\leq\frac{b_1}2$ 
and $\kappa>1$.
Then
\begin{equation}
\label{eqqq}
   \frac{b_1}2\rho^{-\beta}
          \leq \int_{\rho\leq|u|\leq{}\kappa\rho} \psi(u)\,du,\qquad \rho\geq1.
\end{equation}

In place of \eqref{eqsaf10} we now have an {\em equality}  
and therefore 
\begin{equation}
\label{eqsaf20prim}
\begin{aligned}
      \null\|PA_{\lambda}(I-P)\|_{\HS}^2
    &= \lambda^{2d}\,c_{10}
              \int\psi\big(2\lambda{}x\big) \,(\omega_2[\chi_\Omega](2x))^2\,dx\\
      &\geq \lambda^{2d}\,c_{10} \int_{1/\lambda\leq{}|x|\leq1}
                            \psi\big(2\lambda{}x\big) \,(\omega_2[\chi_\Omega](2x))^2\,dx\\
      &\geq \lambda^{2d}\,c_{11} \int_{1/\lambda\leq{}|x|\leq1}
                            \psi\big(2\lambda{}x\big) |x|^{\beta_\Omega}\,dx\\
      &\geq \lambda^{d-\beta_\Omega}\,c_{12} \int_{2\leq{}|u|\leq2\lambda}
                            \psi(u) |u|^{\beta_\Omega}\,du
\end{aligned}
\end{equation}
where we have use the lower estimate in \eqref{eq_ts}
for $\chi_\Omega$. Let $\lambda\geq\kappa$,
set $L=[\log_\kappa\lambda]$ where $[\cdot]$ denotes the integer part
of a number.
Now we split the domain of integration
as a union of concentric domains in the standard way (see e.g.~\cite{BCT})
$$
\begin{aligned}
   \lambda^{d-\beta_\Omega}\,c_{12} \int_{2\leq{}|u|\leq2\lambda}
                            \psi(u) |u|^{\beta_\Omega}\,du
          &\geq\lambda^{d-\beta_\Omega}\,c_{12}
                 \sum_{l=0}^L\int_{2\kappa^l\leq{}|u|\leq\kappa\cdot(2\kappa^l)}
                            \psi(u) |u|^{\beta_\Omega}\,du\\
        &\geq\lambda^{d-\beta_\Omega}\,c_{12}
            \sum_{l=0}^L (2\kappa^l)^{\beta_\Omega}
               \cdot\frac{b_1}2(2\kappa^l)^{-\beta}\\
      &=\lambda^{d-\beta_\Omega}\,c_{13}
            \sum_{l=0}^L \big(\kappa^{\beta_\Omega-\beta}\big)^l
\end{aligned}
$$
where we have used \eqref{eqqq}. Considering the three cases
when $\beta_\Omega$ is smaller than, equal to, and greater than $\beta$
separately, and using the 
fact that $\kappa^{-1}\lambda<\kappa^L\leq\lambda$
together with \eqref{eqsaf20prim} we conclude that 
for the operator $A_\lambda$
\begin{equation}
\label{eq_sharpness}
\|PA_\lambda{}(I-P)\|_{\HS}^2
    \geq c_{14}\cdot\begin{cases} 
             \lambda^{d-\min(\beta,\beta_\Omega)},  &\beta\neq\beta_\Omega\cr
              \lambda^{d-\beta}\log\lambda,     &\beta=\beta_\Omega
       \end{cases}       
\end{equation}
for $\lambda\geq\kappa$.
The sharpness of the estimate \eqref{CH3_s1e2223} is proved.

3. Finally we prove \eqref{CH3_s1e2222}.
Analogously to \eqref{eqsaf10}, \eqref{eqsaf19} we obtain
$$
\begin{aligned}
       \|P{}{}A_{\lambda}\|_{\HS}^2 
      &\leq     \lp\frac{\lambda}{2\pi}\rp^{2d}
             \int\psi(2\lambda{}x) \,dx\,
                 \int \chi_\Omega(x+y)\, 
                     \chi_\Omega(-x+y)\,dy\\ 
      &\leq (2\pi)^{-2d}\,2^{-d}\,\|\psi\|_{1}\,\meas(\Omega)\cdot\lambda^d,
\end{aligned}
$$
for all $\lambda\geq2$.
This finishes the proof of Theorem \ref{CH3_lem.lap.1}.
\subsection{Proof of Theorem \ref{CH3_s0cor1}}
\label{subsec.2.2}
It is important that in both cases
the fact $PA_\lambda\in\HS$ implies that
the operator \eqref{eqfirst}
is trace-class, and also that in this case
the absolute
value of the trace of \eqref{eqfirst} can be estimated 
as follows: there exists a constant $C(\Omega,\sigma)$ such that
one has, for $\lambda\geq2$,
\begin{equation}
\label{eq1thm1}
    \Big|\tr\big[Pf(P{}A_{\lambda}P)P -
           Pf(A_{\lambda})P\big]\Big|
   \leq{}C(\Omega,\sigma)\cdot{}\tilde{C}(f)
    \cdot \|P{}{}A_{\lambda}(I-P{})\|_\HS
\,\|(I-P)A_{\lambda}P\|_\HS.
\end{equation}
In the case of self-adjoint $A_\lambda$,
\eqref{eq1thm1}
 was proved in \cite[Theorem~1.2]{LS}
(note that $\cup_{0\leq t\leq1}\cup_{\lambda\geq2}
\spec A_{\lambda}\subset{\mathbb{R}}$). 
In the case of analytic $f$ the idea of the 
proof goes back to \cite{Wevd}. 
More precisely,
denote $Q=I-P$
and note that for any $m\in\N$, $m\geq2$
we can write
\begin{equation}
\label{eq111}
\begin{aligned}
   PA_\lambda^m P &= PA_\lambda(P+Q)A_\lambda^{m-1}P
                  = P(A_\lambda P)A_\lambda^{m-1}P
                      +(PA_\lambda)Q A_\lambda^{m-1}P\\
      &= P(A_\lambda P)A_\lambda(P+Q)A_\lambda^{m-2}P
                      +(PA_\lambda)Q A_\lambda^{m-1}P\\
     &= P(A_\lambda P)^2 A_\lambda^{m-2}P
                     +(PA_\lambda)^2Q A_\lambda^{m-2}P
                      +(PA_\lambda)Q A_\lambda^{m-1}P=\cdots=\\
     &= P(A_\lambda P)^m 
            + \sum_{j=1}^{m-1} (PA_\lambda )^{m-j}QA_\lambda^j P.
\end{aligned}
\end{equation}
Also for $j\geq2$
\begin{equation}
\label{eq222}
\begin{aligned}
   A_\lambda^j P &= A_\lambda^{j-1}(P+Q)A_\lambda P
                  = A_\lambda^{j-1}PA_\lambda P 
                        + A_\lambda^{j-1}QA_\lambda P \\
                &= A_\lambda^{j-2}(P+Q)A_\lambda PA_\lambda P
                       + A_\lambda^{j-1}QA_\lambda P \\
           &= A_\lambda^{j-2}(PA_\lambda)^2 P
                       + A_\lambda^{j-2}Q(A_\lambda P)^2
                    + A_\lambda^{j-1}QA_\lambda P =\cdots=\\
         &= A_\lambda (PA_\lambda)^{j-1} P
                   + \sum_{k=1}^{j-1} A_\lambda^k Q (A_\lambda P)^{j-k}P.
\end{aligned}
\end{equation}
Substituting \eqref{eq222} in \eqref{eq111} we obtain
\begin{equation}
\label{eqsum2}
\begin{aligned}
  P(&A_\lambda)^mP - (PA_\lambda P)^m 
          = \sum_{j=1}^{m-1}(PA_\lambda)^{m-j}QA_\lambda^jP\\
               &= \sum_{j=1}^{m-1}(PA_\lambda)^{m-j}Q
                             A_\lambda(PA_\lambda)^{j-1}P 
         + \sum_{j=2}^{m-1}(PA_\lambda)^{m-j}Q 
          \sum_{k=1}^{j-1}A_\lambda^kQ(A_\lambda P)^{j-k}P.
\end{aligned}
\end{equation}
There are $\frac{m(m-1)}{2}$ terms on the right-hand side in \eqref{eqsum2}
 each containing both $PA_\lambda{}Q$ and $QA_\lambda{}P$.
Hence
\begin{equation}
\label{eq_nnnew}
\begin{aligned}
  \big\|(PA_{\lambda}P)^m -  P(A_{\lambda})^m{}P\big\|_{\mathfrak{S}_1} 
         \leq\frac{m(m-1)}{2}&\|A_{\lambda}\|^{m-2}\,
        \|P{}{}A_{\lambda}Q\|_{{\mathfrak{S}}_2}
                   \,\|QA_{\lambda}P\|_{{\mathfrak{S}}_2}.
\end{aligned}
\end{equation}
Recall \eqref{CH3_s0esup}
and note that \eqref{eq1thm1}
(even 
with $\big\|Pf(PA_{\lambda}P)P -  Pf(A_{\lambda}){}P\big\|_{\mathfrak{S}_1}$
in the left-hand side)
holds 
for any $f(z)$ 
analytic on a neighborhood of $\{z:|z|\leq{}S(\sigma)\}$.

Note finally that for $f(z)=z^2$ and a self-adjoint $A_\lambda$
$$
        \tr\big[Pf(P{}A_{\lambda}P)P -
           Pf(A_{\lambda})P\big]
             =\|P{}{}{A_{\lambda}}(I-P{})\|_{\HS}^2.
$$
This together with the sharpness
in Theorem \ref{CH3_lem.lap.1}
 implies sharpness in Theorem \ref{CH3_s0cor1}.
\subsection{Proof of Corollary \ref{CH3_s0cor5}}
\label{subsec.2.3}
We need to prove the case (i).
By Theorem \ref{CH3_s0cor1} we only have to analyze
$\tr(Pf(A_\lambda)P)$, $f$ analytic.
We use the functional calculus developed in \cite{W5}.
Denote by $\op{\tau}$ the depending
on the parameter $\lambda\geq2$ operator with integral kernel
$
          (2\pi)^{-d}
              \int e^{i{}\xi\cdot{}(x-y)} \, \tau(x,\xi/\lambda)\,d\xi
$
and let
$
   Q_{\lambda\Gamma} = \op\chi_{\lambda\Gamma}(\xi)$.
Then 
$            A_{\lambda} = (\op \tau)\,Q_{\lambda\Gamma}.
$
By \cite[Lemma D(iii)]{W5}
for any $\eps>0$ there is a constant $\Lambda(\eps)\geq2$
such that for all $\lambda\geq\Lambda(\eps)$
$$
     \|A_{\lambda }\| \leq (1+\eps)\,\|\tau\|_\infty.
$$
Next, by \cite[Lemma 2]{W5},
for any $\epsilon>0$ there exist 
$C(\eps)$ and $\Lambda(\eps)\geq2$ 
such that for all $m\in{}\N$
and $\lambda\geq\Lambda(\eps)$
\begin{equation}
\label{eqWaux}
\begin{aligned}
      \| P\,(A_\lambda)^m\,P
         - P
     \,(\op \tau^m)\,Q_{\lambda\Gamma}\,P\|_{{\mathfrak{S}}_1}
              \leq C(\eps)\,(1+\eps)^m\,\|\tau\|_\infty^m\cdot\lambda^{d-1}.
\end{aligned}
\end{equation}
Clearly
$
        \tr (P\,(\op \tau^m)\,Q_{\lambda\Gamma}\,P) = 
               \big(\frac{\lambda}{2\pi}\big)^d\int_\Gamma\int_\Omega 
                      (\tau(x,\xi))^m\,dx\,d\xi.
$
By the argument \cite[p.~184]{W5}, as $\lambda\ra\infty$
$$
\begin{aligned}
     \tr (P\,(f(\op \tau))\,Q_{\lambda\Gamma}\,P)
       = \bigg(\frac{\lambda}{2\pi}\bigg)^d\int_\Gamma\int_\Omega 
              f(\tau(x,\xi))\,dx\,d\xi
         +O(\lambda^{d-1}).
\end{aligned}
$$
Now for any $f(z) = \sum_{m=1}^\infty{}c_m\,z^m$ 
analytic on a neighborhood of $\{z:|z|\leq\|\tau\|_\infty\}$
 we have $c_m=O((1+\delta)^{-m}\|\tau\|_\infty^m)$ for some $\delta>0$.
Taking $\epsilon<\delta$ in \eqref{eqWaux} 
finishes the proof of Corollary \ref{CH3_s0cor5}.
\subsection{Proof of Lemma~\ref{lemex}}
\label{subsec.2.4}
We start with $d=1$.
Here we write $\beta=\beta_\Omega$ for brevity
and denote the Lebesgue measure in $\R$ by $\meas_1$.
The case $\beta=1$ is trivial, any finite interval with 
non-empty interior would do.
Let us therefore assume $0<\beta<1$. 
We construct a set $\Omega$ with the required properties
as a (finite union of sets each of which is a) countable 
union of closed intervals
obtained by a process similar to the construction of a Cantor set.
The difference is that we do not remove but rather add intervals.
Let $0<\alpha<\infty$, $I_\alpha=[0,\alpha]$ and $0<q<1$.
We explain how $q$ is related to the given $\beta$ later.
We construct the set $\Omega_\alpha$ as follows.
We start with an empty set and at 
the $0$th step we add the middle $q$th part of $I_\alpha$ 
to $\Omega_\alpha$.
Each of the remaining intervals has length $\alpha{}Q$,
where we denote $Q=\frac{1-q}2$.
Note that $0<Q<\frac12$.
Then we take the middle $q$th part of
each of the two remaining intervals and do {\em not\ }add it
to the set $\Omega_\alpha$. 
Each of the four remaining intervals has length
$\alpha{}Q^2$. 
Now we take the $q$th middle
part of each of the four intervals and add it 
to $\Omega_\alpha$. 
This completes
the $1$st step of the construction of $\Omega_\alpha$. 
We continue in this manner.
The set $\Omega_\alpha$ we obtain
has the following properties:

1. For each $k\in\N$, it contains $2^{2k}$ intervals
of length $\alpha{}qQ^{2k}$;

2. Half of these $2^{2k}$ intervals will have an interval of length
$\alpha{}qQ^{2k-1}$, 
at a distance of $\alpha{}Q^{2k+1}$ to the right, 
which does {\em not\ }contain points from $\Omega_\alpha$.

\noindent
For any $k\in\N$ introduce notation 
\begin{equation}\label{eqnot}
\begin{aligned}
   a_k^{(\alpha)}&= \alpha(qQ^{2k}+Q^{2k+1}) 
       = \alpha(q+Q)Q^{2k}\\
   b_k^{(\alpha)}&= \alpha(Q^{2k+1}+qQ^{2k-1})
       = \alpha(qQ^{-1}+Q)Q^{2k}
\end{aligned}
\end{equation}
and consider the
shift of the set $\Omega_\alpha$ to the right by 
\begin{equation}
\label{eqex}
        h\in\big[  a_k^{(\alpha)},\, b_k^{(\alpha)} \big],\qquad k\in\N,
\end{equation}
units.
Recall that $Q\in(0,1/2)$ and 
set $\beta=1-\frac{\log2}{\log1/Q}\in(0,1)$.
(Actually we start with a given $\beta\in(0,1)$ and after that
define $Q$ as above and $q=1-2Q\in(0,1)$.)
By the properties 1 and 2 above for $h$ as in \eqref{eqex}
\begin{equation}
\label{eqstar}
  \meas_1\big((\Omega_\alpha-h)\setminus\Omega_\alpha\big)
       \geq \frac12\cdot2^{2k}\cdot\alpha{}qQ^{2k}
        =\frac12\alpha{}q \cdot\big(Q^{2k}\big)^{1-\frac{\log2}{\log1/Q}}.
\end{equation}
We note that in \eqref{eqex},
$
       c_1  Q^{2k}\leq h\leq c_2 Q^{2k}
$,
where $c_1=\alpha(q+Q)$ and $c_2=\alpha(qQ^{-1}+Q)$ 
do not depend on $k$. 
In particular $Q^{2k}\geq\frac1{\alpha(qQ^{-1}+Q)}\,h$
for $h$ as in \eqref{eqex}, and hence \eqref{eqstar} implies
\begin{equation}
\label{eqex1}
      \meas_1\big((\Omega_\alpha-h)\setminus\Omega_\alpha\big)
      \geq \frac{\alpha{}q}{2(\alpha(qQ^{-1}+Q))^{\beta}}\cdot h^\beta
\end{equation}
for $h$ as in \eqref{eqex}, and therefore 
for $h\in\bigcup_{k\in\N}[a_k^{(\alpha)},\, b_k^{(\alpha)} \big]$. 
Using $q=1-2Q$ we find that
$$
    \frac{b_{k+1}^{(\alpha)}}{a_k^{(\alpha)}} 
      = \frac{qQ+Q^3}{q+Q} = Q-Q^2\in(0,1/4),\qquad q\in(0,1),
$$
and hence the set
$\bigcup_{k\in\N}[a_k^{(\alpha)},\, b_k^{(\alpha)} \big]$ has gaps.
We need \eqref{eqex1} to hold for {\em all} small $h$. 
In order to fulfill this condition, we consider a finite union of the scaled sets
$\Omega_\alpha$ for appropriate $\alpha$. 
More precisely, the first set in the union is the set $\Omega_1$
corresponding to $\alpha=1$. 
Then \eqref{eqex1} holds for $h\in\bigcup_{k\in\N}[a_k^{(1)},b_k^{(1)}]$.
Note that the ratio
\begin{equation}\la{eqgamma}
    \gamma= \frac{b_k^{(\alpha)}}{a_k^{(\alpha)}}=\frac{q+Q^2}{qQ+Q^2} >1
\end{equation}
is independent of $k$ and $\alpha$. 
Now for this $\gamma$ we construct the set $\Omega_\gamma$. 
For the set $\Omega_\gamma$,
\eqref{eqex1} holds (with a different constant) for
$
     h\in\bigcup_{k\in\N}[\gamma{}a_k^{(1)},\gamma{}b_k^{(1)}]
        =\bigcup_{k\in\N}[b_k^{(1)},\gamma{}b_k^{(1)}]
$
because $\gamma{}a_k^{(\alpha)}=b_k^{(\alpha)}$ for all $k\in\N$
and $\alpha>0$, in particular for $\alpha=\gamma$.
Now for the union $\Omega_1\cup(\Omega_\gamma-2)$,
 \eqref{eqex1} holds (with a different constant)
for 
$$
   h\in\bigg(\bigcup_{k\in\N}[a_k^{(1)},b_k^{(1)}]\bigg)\cup
       \bigg(\bigcup_{k\in\N}[b_k^{(1)},\gamma{}b_k^{(1)}]\bigg)
        =\bigcup_{k\in\N}[a_k^{(1)},\gamma{}b_k^{(1)}].
$$
Now we construct
the scaled sets $\Omega_{\gamma^2},\cdots,\Omega_{\gamma^N}$,
where the (finite) number $N\in\N$ is determined by the (independent of $k$ 
and $\alpha$)
condition $\gamma^N b_{k}^{(\alpha)}\geq a_{k-1}^{(\alpha)}$ or
$\gamma^N\geq\frac{q+Q}{qQ+Q^3}$
(recall that $\gamma>1$ by \eqref{eqgamma}).
Set finally
$$
     \Omega= \Omega_1\cup(\Omega_\gamma-2)
    \cup(\Omega_{\gamma^2}-3-\gamma)\cup\cdots
    \cup\bigg(\Omega_{\gamma^N}-\sum_{j=1}^N(1+\gamma^{j-1})\bigg).
$$
Since $N$ is finite, \eqref{eqex1} for the set $\Omega$ 
holds (with a different constant) for $h$ in
$$
\begin{aligned}
     \bigcup_{k=2}^\infty&\bigg([a_k^{(1)},b_k^{(1)}]
             \cup[b_k^{(1)},\gamma{}b_k^{(1)}]
              \cup[\gamma{}b_k^{(1)},\gamma^2b_k^{(1)}]
              \cup\cdots\cup[\gamma^{N-1}b_k^{(1)},\gamma^Nb_k^{(1)}]\bigg)\\
      &\supset\bigcup_{k=2}^\infty\big[a_k^{(1)},a_{k-1}^{(1)}\big]
        =\big(0,a_1^{(1)}\big].
\end{aligned}
$$
Therefore \eqref{eqex1} holds for the set $\Omega$
with some constant for  {\em all\ } $h\in(0,a_1^{(1)}]$.

Now we explain why the 
estimate opposite to \eqref{eqex1} holds for the same $\beta$.
Again since $N$ is finite
it suffices to consider one set $\Omega_\alpha$
with an arbitrary $\alpha>0$.
Choose first any $k=2,3,\cdots$ 
and $h\in[a_k^{(\alpha)},a_{k-1}^{(\alpha)}]$. 
We estimate the contribution of the sets up to
generation $k$ to $\meas_1((\Omega_\alpha-h)\setminus\Omega_\alpha)$
by $a_{k-1}^{(\alpha)}$ times their total number.
We estimate the conribution of the sets from generation $k+1$
and onwards by their total length. As a result we obtain
the following estimate
\begin{equation}
\label{labex3}
\begin{aligned}
   \meas_1\big((\Omega_\alpha-h)\setminus\Omega_\alpha\big)
          &\leq a_{k-1}^{(\alpha)} \sum_{l=0}^k2^{2l}
              + \sum_{l=k+1}^\infty \alpha{}qQ^{2l} \cdot2^{2l} \\
         &\leq c\cdot{}Q^{2k}\cdot2^{2k}
         =c\cdot \big(Q^{2k}\big)^{1-\frac{\log2}{\log1/Q}} 
       \leq\tilde{c}\cdot h^\beta
\end{aligned}
\end{equation}
for some $c,\tilde{c}$ independent of $k$ (recall \eqref{eqnot}).
Since the right-hand side in \eqref{labex3}
does not depend on $k$, we conclude that
 \eqref{labex3} holds 
for $h\in\bigcup_{k=2}^\infty[a_k^{(\alpha)},a_{k-1}^{(\alpha)}]
=(0,a_1^{(\alpha)}]$.
This completes the proof for $d=1$.

In the case $d\geq2$ for a given $\beta\in(0,1)$,
let $\Omega$ be the constructed above set
and consider the direct product $\Omega^d$.
Let $\meas_d$ denote the Lebesgue measure in $\R^d$.
Recall that $(\Omega^d)_h=\Omega^d-h$, $h\in\R^d$.
Since for an arbitrary $h=(h_1,\cdots,h_d)$ we can write $\Omega^d$
as a disjoint union
$
         \Omega^d = (\Omega^d \cap (\Omega^d)_h )
            \cup (\Omega^d \setminus (\Omega^d)_h )
$
we have 
\be\la{eq1}
        \meas_d\big(\Omega^d \setminus (\Omega^d)_h \big) =
                  \meas_d\big(\Omega^d\big) 
         -  \meas_d\big( \Omega^d \cap (\Omega^d)_h \big).
\ee
Noting that $\Omega^d \cap (\Omega^d)_h = 
(\Omega\cap\Omega_{h_1})\times\cdots \times(\Omega\cap\Omega_{h_d})$
and using \eqref{eq1} for each of the factors we find
\be\la{eq2}
\ba
      \meas_d\big(\Omega^d \cap (\Omega^d)_h\big) &=
         \prod_{j=1}^d \meas_1\big(\Omega\cap\Omega_{h_j}\big)\\ 
             &= \prod_{j=1}^d \Big( \meas_1(\Omega) 
                 - \meas_1\big(\Omega\setminus\Omega_{h_j}\big) \Big).
\ea
\ee
By the construction of the set $\Omega$, there are $c_1,c_2>0$
so that for $|h_j|$ small,
$$
      c_1|h_j|^\beta\leq
              \meas_1\big(\Omega\setminus\Omega_{h_j}\big) \leq 
                  c_2 |h_j|^\beta
$$
which together with \eqref{eq2} implies
\be\la{eq2.5}
      \prod_{j=1}^d \big(\meas_1(\Omega) - c_2 |h_j|^\beta \big)
           \leq \meas_d\big(\Omega^d \cap (\Omega^d)_h\big)
           \leq \prod_{j=1}^d \big(\meas_1(\Omega) - c_1 |h_j|^\beta \big).
\ee
Recall that $|h|=(\sum_{j=1}^d|h_j|^2)^{1/2}$, $h\in\R^d$
(we hope that denoting the Euclidean length of a vector
and the absolute value of a number does not lead to confusion below).
Note that 
\be\la{eq3}
    \prod_{j=1}^d \Big(\meas_1(\Omega) - c_1 |h_j|^\beta \big)
       = \big(\meas_1(\Omega)\big)^d   
             - c_1 \sum_{j=1}^d |h_j|^\beta + E(h_1,\cdots,h_d)   
\ee
where the function $E$ has the property that it is a finite sum
of terms each of which contains a factor
 $|h_{k_1}|^\beta\cdots|h_{k_p}|^\beta $
for some $1\leq k_1<\cdots<k_p\leq d$ and $p\geq2$.
Note that each such factor is $\leq O\big(|h|^{2\beta}\big)$ for small $|h|$.
(Indeed, we can write 
$
       |h_1|^\beta |h_2|^\beta = |h|^{2\beta} 
       \cdot (|h_1|/|h|)^\beta
          \cdot(|h_2|/|h|)^\beta \leq |h|^{2\beta}
$
and for all other possible factors 
again use the estimate $|h_j|\leq |h|$.)
Substituting \eqref{eq3} in \eqref{eq2.5} and
returning to \eqref{eq1}, 
 we find 
\be\la{eq4}
       c_1 \sum_{j=1}^d |h_j|^\beta  + O(|h|^{2\beta})
       \leq \meas_d\big(\Omega^d \setminus (\Omega^d)_h\big)
            \leq  c_2 \sum_{j=1}^d |h_j|^\beta  + O(|h|^{2\beta})
\ee
for small enough $|h|$ (we have used 
that $\meas_d(\Omega^d)=(\meas_1(\Omega))^d$).
We note finally that for any fixed $\beta\in(0,1)$ there exist $C_1,C_2>0$
so that for all $h$
\be\la{eq5}
     C_1 |h|^\beta \leq \sum_{j=1}^d |h_j|^\beta \leq C_2 |h|^\beta.
\ee
(To see this, set $\tilde{h}=h/|h|$ and note that we just have to 
prove that for some $\tilde{C}_{1},\tilde{C}_2>0$,
$
     \tilde{C}_1 \leq \sum_{j=1}^d |\tilde{h}_j|^\beta 
              \leq \tilde{C}_2$
if only $\sum_{j=1}^d|\tilde{h}_j|^2=1$.
But this holds since the minimum and the maximum are attained 
at some points on the unit sphere.)
Combining \eqref{eq4} and \eqref{eq5}
completes the proof of Lemma~\ref{lemex}.
\par
{\bf Acknowledgments.}
The author would like to thank the following colleagues
for useful discussions: A.~Laptev for the suggestion
 to work on the Widom
conjecture,  Yu.~Safarov for a simpler and more general proof of
Theorem~\ref{CH3_lem.lap.1} than the original one in \cite{Gi},
S.~Smirnov
for the idea of the construction
of the set in Lemma~\ref{lemex},
and I.~Klich, who represented the EE in a form amenable
to the analysis related to the Widom conjecture (see \cite{klich}).
The work in this paper was supported in part by the Royal
Institute of Technology (KTH) in Stockholm,
the Swedish Foundation for International Cooperation
in Research and Higher Education (STINT) grant Dnr.~PD2001--128,
the NSF grant INT--0204308 U.S.--Sweden Collaboratiove Workshop
on PDE's and Spectral Theory, and the NSF grant DMS--0556049.

\end{document}